\documentclass[11pt]{amsart}
\title{Towards commutator theory for relations. II}
\keywords{Commutator, congruence, tolerance, relation, neutral}
\subjclass[2000]{Primary 08A30; Secondary 08B10}
\author{Paolo Lipparini}
\address{Dipartimento di Matematica, Viale della Ricerca Scientifica,
II Universit\`a di Roma (Tor Verguta),
 ROME 
ITALY}
\thanks{The author has received support from MPI and GNSAGA} 

\email{lipparin@axp.mat.uniroma2.it}
\urladdr{http://www.mat.uniroma2.it/\textasciitilde lipparin}
\newtheorem{Theorem}{Theorem}
\newtheorem{proposition}[Theorem]{Proposition}

\newtheorem{thm}[Theorem]{Theorem}
\newtheorem{theorem}[Theorem]{Theorem}

\newtheorem{lemma}[Theorem]{Lemma}

\theoremstyle{definition}
\newtheorem{problem}[Theorem]{Problem}
\newtheorem{remark}[Theorem]{Remark}

\newcommand{\alg}{\mathbf} 
\def\v{\mathcal V}  

\newcommand{\adma}{\mathrm{Adm({\alg A})}} 
\newcommand{\admx}{\mathrm{Adm({\alg X})}} 

\newcommand{\cona}{\mathrm{Con({\alg A})}} 

\begin{document}

\begin{abstract} 
We find conditions equivalent to some commutator 
identities considered in Part I
\end{abstract} 

\maketitle
\bigskip 

See Part I \cite{L}  for notations.

\begin{theorem}\label{x32} 
For every variety $\v$, the following are equivalent:
\[ 
\tag{i}
R \subseteq [R, R|1]
\]
for every algebra ${\alg A} \in \v$ and
for every 
reflexive compatible relation $R$ on ${\alg A} $.
\[ 
\tag{ia}
R^* \subseteq [R, R|1]
\]
for every algebra ${\alg A} \in \v$ and
for every 
reflexive compatible relation $R$ on ${\alg A} $.
\[
\tag{ii}
R \cap T \subseteq [R, T |1]
\]
for every algebra ${\alg A} \in \v$ and
for all reflexive compatible relations $T,R$ on ${\alg A} $.
\[
\tag{iii}
 (R_1 \circ R_2)\cap T \subseteq 
\big(
T \cap (R_2^- \circ (T \cap (R_1^- \circ R_1)) \circ R_2)
\big)^*
\]
for every algebra ${\alg A} \in \v$ and
for all
reflexive compatible relations $R_1$, $R_2$, $T$ on ${\alg A} $.
\[
\tag{iv}
 R_1 \cap (T \circ R_2) \subseteq
\big(
T \cap (R_2 \circ (T \cap (R_1^- \circ R_1)) \circ R_2^-)
\big)^*
\circ R_2
\] 
for every algebra ${\alg A} \in \v$ and
for all
reflexive compatible relations $R_1$, $R_2$, $T$ on ${\alg A} $.
\[
\tag{v}
\beta  \cap (T \circ S) \subseteq
\big(
T \cap (S \circ (T \cap \beta ) \circ S)
\big)^*
\circ S
\] 
for every algebra ${\alg A} \in \v$ and
for every congruence $ \beta $, tolerance $S$
and reflexive compatible relation $T$ on ${\alg A} $.
\[
\tag{vi}
\beta  \cap (T \circ \gamma) \subseteq
(\gamma \circ (T \cap \beta ))^*
\] 
for every algebra ${\alg A} \in \v$ and
for all congruences $ \beta, \gamma $
and every reflexive compatible relation $T$ on ${\alg A} $.

\item[ (vii)] For some $n$, $\v$ has $4$-ary terms $f_i$, $i=0,\dots,n$ such that 
\[ 
 x=f_0(x,y,z,x), \qquad f_n(x,y,z,y)=z, \qquad 
\]
\[ 
 f_i(x,y,x,x)= f_i(x,y,x,y),  \qquad \text{for } i=0,\dots,n,
 \]
\[ 
f_{i-1}(x,y,y,y)= f_i(x,y,y,x), \qquad \text{for } i=1,\dots,n
\]
\begin{multline*} 
\tag{viii}
 (R_1 \circ R_2 \circ \dots \circ R_{n-1} \circ R_n)\cap T \subseteq 
\\
\big(
T \cap (R_n^- \circ (T \cap (R_{n-1}^- \circ \dots
(T \cap (R_2^- \circ (T \cap (R_1^- \circ R_1)) \circ R_2))
\dots \circ R_{n-1})) \circ R_{n})
\big)^*
\end{multline*}
for every algebra ${\alg A} \in \v$ and
for all
reflexive compatible relations $R_1$,\dots, $R_n$ and $T$ on ${\alg A} $.
\end{theorem} 

\begin{proof}
In view of Part I, Theorem 3, if we prove
(vi) $ \Rightarrow $ (vii) $\Rightarrow $ (i)
then (i)-(vii) are all equivalent.

(vi) $ \Rightarrow $ (vii) is a slight variation on an argument
from \cite{KK}. 
Consider the free algebra ${\alg F} $ in $\v$
generated by the three elements $x$, $y$, $z$,
and let $\beta=Cg(x,z)$, $\gamma=Cg(y,z)$, and
$T$ be the smallest {\em admissible reflexive relation}
containing $(x,y)$.

Thus, $(x,z)\in \beta \cap(T \circ \gamma )$ and, by (vi),
$(x,z)\in ((T \cap \beta ) \circ \gamma)^*$.
This means that there is an integer $n$
and there are terms $t_0(x,y,z)$,\dots, 
$t _{2n+1}  (x,y,z)$ such that 
\[ 
 x=t_0(x,y,z), \qquad t _{2n+1} (x,y,z)=z, \qquad 
\]
\[ 
 (t _{2i}  (x,y,z), t _{2i+1}  (x,y,z)) \in T \cap \beta, 
 \qquad \text{for } i=0,\dots,n, \text{ and} 
 \]
\[ 
 (t _{2i-1}  (x,y,z), t _{2i}  (x,y,z)) \in \gamma, 
 \qquad \text{for } i=1,\dots,n
 \]
Notice that $ T= \{(f(x,y,z,x), f(x,y,z,y))| f \text{ a term of } 
{\alg F} \} $, since the right-hand relation is  reflexive,
admissible, and contains $(x,y)$; moreover, every other reflexive admissible
relation containing $(x,y)$ contains all pairs of the form
$(f(x,y,z,x), f(x,y,z,y))$.

Hence there are terms $f_i(x,y,z,w)$, $i=0,\dots, n$
such that 
\[ 
 x=t_0(x,y,z), \qquad t _{2n+1} (x,y,z)=z, \qquad 
\]
\[ 
t _{2i}  (x,y,z)=f_i(x,y,z,x), \quad f_i(x,y,z,y)= t _{2i+1}  (x,y,z),  
 \quad \text{for } i=0,\dots,n, 
 \]
\[ 
t _{2i}  (x,y,x)= t _{2i+1}  (x,y,x),  
 \qquad \text{for } i=0,\dots,n, \text{ and} 
 \]
\[ 
 t _{2i-1}  (x,y,y)= t _{2i}  (x,y,y), 
 \qquad \text{for } i=1,\dots,n
 \]
By using the identities in the second line, we can express
all the other identities by means of the $f_i$'s, thus getting
the desired relations. 

 (vii) $\Rightarrow $ (i).
Let ${\alg A} $ be an algebra $ \in \v$ and
 let $R$ be a reflexive compatible relation on ${\alg A} $.
If $x,y \in {\alg A} $ and $x R y$ then for every $i$ the matrix
\[
\begin{vmatrix}   
f_i(x,y,x,x) & f_i(x,y,x,y) \cr
f_i(x,y,y,x) & f_i(x,y,y,y) 
\end{vmatrix}  
\] 
belongs to $M(R,R)$, hence for every $i$
$ f_i(x,y,y,x) [R,R|1]  f_i(x,y,y,y) $,
since $ f_i(x,y,x,x) = f_i(x,y,x,y) $.

Thus,     $x= f_0(x,y,y,x)
 [R,R|1]  
f_0(x,y,y,y)
= f_1(x,y,y,x) 
[R,R|1] 
f_1(x,y,
$ $
y,y)
=
f_2(x,y,y,x) \dots  f _{n-1} (x,y,y,y)=
f_n(x,y,y,x) [R,R|1]  f_n(x,y,y,y)=y$.
This implies that for every $x,y \in {\alg A} $,
if $x R y $  then
$x [R,R|1] y $, that is, 
$ R \subseteq [R,R|1]$. 

(viii) $\Rightarrow $ (iii) is trivial: just let
$R_n$, $R_{n-1}$,\dots, $R_3$ be equal to  the identity
relation.

(ii) $\Rightarrow $ (viii) is immediate from the next Lemma.
\end{proof}

\begin{lemma}\label{x1c}
For $R_1$,\dots, $R_n$, $S$, $T$ reflexive compatible relations  on some algebra,
the following hold:
\begin{multline*} 
\tag{i}
 [R_1 \circ R_2 \circ \dots \circ R_{n-1} \circ R_n, T|1] \subseteq 
\\
\big(
T \cap (R_n^- \circ (T \cap (R_{n-1}^- \circ \dots
(T \cap (R_2^- \circ (T \cap (R_1^- \circ R_1)) \circ R_2))
\dots \circ R_{n-1})) \circ R_{n})
\big)^*
\end{multline*}
\begin{multline*} 
\tag{ii}
K(R_1 \circ R_2 \circ \dots \circ R_{n-1} \circ R_n, T;S) \subseteq 
\\
K(R_n,T; K(R _{n-1},T; \dots K(R_2,T;K(R_1,T;S))\dots))
\end{multline*}
Notice that (ii) can be restated as follows: define 
$K_1= K(R_1,T;S)$, and 
$K_{i+1}= K(R_{i+1},T;K_i)$. Then 
\[   
K(R_1 \circ R_2 \circ \dots \circ R_{n-1} \circ R_n, T;S) \subseteq K_n
\] 
\end{lemma}
\begin{proof} Compare Lemmata 1(i) and 6(i) in Part I.
\end{proof}

\begin{theorem}\label{x22} 
For every variety $\v$, the following are equivalent:
\[ 
\tag{i}
R \subseteq [R, R^\circ|1]
\]
for every algebra ${\alg A} \in \v$ and
for every 
reflexive compatible relation $R$
 on ${\alg A} $.
\[ 
\tag{ia}
R^* \subseteq [R, R^\circ|1]
\]
for every algebra ${\alg A} \in \v$ and
for every 
reflexive compatible relation $R$ on ${\alg A} $.
\[ 
\tag{ib}
R^-\subseteq [R, R^\circ|1]
\]
for every algebra ${\alg A} \in \v$ and
for every 
reflexive compatible relation $R$ on ${\alg A} $.
\[ 
\tag{ic}
R^\circ \subseteq [R, R^\circ|1]
\]
for every algebra ${\alg A} \in \v$ and
for every 
reflexive compatible relation $R$ on ${\alg A} $.
\[ 
\tag{id}
Cg(R) = [R, R^\circ|1]
\]
for every algebra ${\alg A} \in \v$ and
for every 
reflexive compatible relation $R$ on ${\alg A} $.
\[
\tag{ii}
R \cap T \subseteq [R, T |1]
\]
for every algebra ${\alg A} \in \v$ and
for every tolerance $T$ and every
reflexive compatible relation $R$ on ${\alg A} $.
\[
\tag{iii}
 (R_1 \circ R_2)\cap T \subseteq 
\big(
T \cap (R_2^- \circ (T \cap (R_1^- \circ R_1)) \circ R_2)
\big)^*
\]
for every algebra ${\alg A} \in \v$ and
for every tolerance $T$ and all
reflexive compatible relations $R_1$, $R_2$ on ${\alg A} $.
\[
\tag{iv}
 R_1 \cap (T \circ R_2) \subseteq
\big(
T \cap (R_2 \circ (T \cap (R_1^- \circ R_1)) \circ R_2^-)
\big)^*
\circ R_2
\] 
for every algebra ${\alg A} \in \v$ and
for every tolerance $T$ and all
reflexive compatible relations $R_1$, $R_2$ on ${\alg A} $.
\[
\tag{v}
\beta  \cap (T \circ S) \subseteq
\big(
T \cap (S \circ (T \cap \beta ) \circ S)
\big)^*
\circ S
\] 
for every algebra ${\alg A} \in \v$ and
for every congruence $ \beta $ and tolerances $T, S$ on ${\alg A} $.
\[
\tag{vi}
\beta  \cap (T \circ \gamma) \subseteq
\gamma \vee (T \cap \beta )^*
\] 
for every algebra ${\alg A} \in \v$ and
for all congruences $ \beta, \gamma $ and tolerance $T$ on ${\alg A} $.
\item[ (vii)] For some $n$, $\v$ has $5$-ary terms $f_i$, $i=0,\dots,n$ such that 
\[ 
 x=f_0(x,y,z,x,y), \qquad f_n(x,y,z,y,x)=z, \qquad 
\]
\[ 
 f_i(x,y,x,x,y)= f_i(x,y,x,y,x),  \qquad \text{for } i=0,\dots,n,
 \]
\[ 
f_{i-1}(x,y,y,y,x)= f_i(x,y,y,x,y), \qquad \text{for } i=1,\dots,n
\]
\begin{multline*}
\tag{viii}
 (R_1 \circ R_2 \circ \dots \circ R_{n-1} \circ R_n)\cap T \subseteq 
\\
\big(
T \cap (R_n^- \circ (T \cap (R_{n-1}^- \circ \dots
(T \cap (R_2^- \circ (T \cap (R_1^- \circ R_1)) \circ R_2))
\dots \circ R_{n-1})) \circ R_{n})
\big)^*
\end{multline*}
for every algebra ${\alg A} \in \v$ and
for all
reflexive compatible relations $R_1$,\dots, $R_n$, and tolerance 
$T$ on ${\alg A} $.
\end{theorem} 

\begin{proof}
The proof is similar to the proof of Theorem \ref{x32}, using
Part I, Theorem 2.
(vi) $\Leftrightarrow $ (vii) is due to \cite{KK}.
\end{proof}

\begin{remark}\label{var} 
We could have defined
$[ R, S | 1]$
to be the {\em smallest congruence} containing
the set   
\[\left\{ 
(z,w) | 
\begin{vmatrix}   
x & x \cr z & w 
\end{vmatrix}  
\in M( R , S ) 
\right\} 
\] 
(rather than its {\em transitive closure}).
With this definition, Theorem \ref{x32} 
still holds, provided $(\ -\ )^*$ is replaced everywhere by $Cg(\ -\ )$,
and condition (vii) is appropriately modified.  
Notice that 
$[R,S|1]^-= [R,S^-|1]$,
$K(R,S;T)^-= K(R,S^-;T^-)$,
and $Cg( [R,S|1] )=
\big(
[R,S|1] \circ [R,S|1]^-
\big)^*=
\big(
K(R,S;0) \circ K(R,S^-;0)
\big)^*
$.
Notice also that, say, 
$Cg( \gamma \circ (T \cap \beta ))=
\big(
\gamma \circ (T \cap \beta )  \circ (T^- \cap \beta )
\big)^*
$, if $\gamma $ and $\beta $ are congruences and
$T$ is reflexive and admissible.  
\end{remark} 

\begin{theorem}\label{x32var} 
We still get equivalent conditions if in the
statement of Theorem 
\ref{x32}:
we replace $ [-,-|1] $ by $ Cg([-,-|1]) $ in conditions (i), (ia), (ii), 
{\em and} 
we replace $ (-)^* $ by $ Cg(-) $ in conditions (iii)-(vi), (viii), 
{\em and} 
we modify condition (vii) to:
\item[ (vii)$'$] For some even $n$, $\v$ has $4$-ary terms $f_i$, $i=0,\dots,n$ such that 
\[ 
 x=f_0(x,y,z,x), \qquad f_n(x,y,z,y)=z, \qquad 
\]
\[ 
 f_i(x,y,x,x)= f_i(x,y,x,y),  \qquad \text{for } i \text{ even, }  0 \leq i \leq n,
 \]
\[ 
f_{i-1}(x,y,y,y)= f_i(x,y,y,y), \qquad \text{for } i \text{ odd, }  1 \leq i \leq n,
\]
\[ 
 f_i(x,y,x,y)= f_i(x,y,x,x),  \qquad \text{for } i \text{ odd, }  0 \leq i \leq n,
 \]
\[ 
f_{i-1}(x,y,y,x)= f_i(x,y,y,x), \qquad \text{for } i \text{ even, }  1 \leq i \leq n
\]
\end{theorem} 

\begin{problem} \label{chr}  
Classify varieties according to the commutator identities for relations they
satisfy. For example, the following properties can be taken into account:

(i) $R \subseteq [R, R^\circ|1]$,

(ii) $R \subseteq [R, R|1]$,

(iii) $R^- \subseteq [R, R|1]$,

(iv) $R \subseteq Cg([R, R|1])$,

(v) $R \subseteq [R, R]$,

(vi) there is a (weak) difference term with respect to $[R, R^\circ|1]$,

(vii) there is a (weak) difference term with respect to $[R, R|1]$,

(viii) there is a (weak) difference term with respect to $Cg([R, R|1])$,

(ix) there is a (weak) difference term with respect to $[R, R]$.

We have partial results suggesting that a classification as above 
agrees with D. Hobby and R. McKenzie's classification of locally finite 
varieties, as extended in most cases to arbitrary varieties by K. Kearnes 
and others.
  \end{problem}

\def\cprime{$'$} \def\cprime{$'$}
\providecommand{\bysame}{\leavevmode\hbox to3em{\hrulefill}\thinspace}
\providecommand{\MR}{\relax\ifhmode\unskip\space\fi MR }
\providecommand{\MRhref}[2]{%
  \href{http://www.ams.org/mathscinet-getitem?mr=#1}{#2}
}
\providecommand{\href}[2]{#2}

\end{document}